\documentclass[12pt]{amsart}
\usepackage{amssymb}
\usepackage{amsfonts}
\usepackage{amsmath}

\allowdisplaybreaks[2]

\newcommand{\bb}{\begin{equation}}
\newcommand{\ee}{\end{equation}}
\newcommand{\bq}{\begin{eqnarray}}
\newcommand{\eq}{\end{eqnarray}}
\newcommand{\bqn}{\begin{eqnarray*}}
\newcommand{\eqn}{\end{eqnarray*}}

\newtheorem{theorem}{Theorem}[section]

\theoremstyle{definition}

\newtheorem{remark}[theorem]{Remark}
\newtheorem{definition}[theorem]{Definition}

\theoremstyle{remark}

\renewcommand{\eqref}[1]{(\ref{#1})}

\renewcommand{\bigskip}{\vspace{0.2cm}}

\renewcommand{\div}{{\rm div}}

\begin{document}

\title{A mathematical clue to  the separation  phenomena on the
two-dimensional Navier-Stokes equation
}
\author{}
\address{}
\email{}
\author{Tsuyoshi Yoneda}
\address{Department of Mathematics, Hokkaido University, Sapporo 060-0810, Japan}
\email{yoneda@math.sci.hokudai.ac.jp}

\date{\today}
\maketitle
\begin{abstract}
In general,
 before separating from a boundary, the flow moves toward reverse direction near the
boundary against the laminar flow direction.
Here in this paper,  a clue to  such reverse flow phenomena (in the mathematical sense) is observed. More precisely, the non-stationary two-dimensional Navier-Stokes equation with an initial datum having a parallel laminar  flow (we define it rigorously in the paper) is considered.
We show that the  direction of the material differentiation is opposite to the initial flow direction
and
 effect of the material differentiation (inducing the reverse flow) becomes bigger
when
the
curvature of the boundary becomes bigger.
We also show that the parallel laminar flow cannot be a stationary Navier-Stokes flow near a portion of the boundary with nonzero curvature.
\end{abstract}
\par

\maketitle

\section{Introduction}
\noindent
Uchida, Sugitani and Ohya  \cite{USO} studied a non-stratified airflow past a two-dimensional ridge in a uniform flow.
Airflows around the ridge include an unsteady vortex shedding. Their study  focused on airflow characteristics in a wake region.
In general,
 before separating from a boundary, the flow moves toward reverse direction near the
boundary against the laminar flow direction (see \cite{JH,USO} for example).
There are several results related to the wake region in pure mathematics. Using the Oseen system is one of the mathematical approach to analyze the wake region.
For the detailed discussion of the Oseen system, we refer the reader to \cite{G}.
 In  a convex obstacle case,  the character of
the system is elliptic in front of the obstacle. To the contrary,  its character changes into
parabolic type (wake region) behind the obstacle (see \cite{K} for example).
Maekawa \cite{M} considered the two-dimensional Navier-Stokes equations in a half plane under the
no-slip boundary condition. He established a solution formula for the vorticity equations
and got a  sufficient condition on the initial data for the vorticity to blow up to the inviscid limit (see also \cite{G}).
 His observation suggests that if the Reynolds number is high,  the boundary layer immediately appears and the high
vorticity creation occurs near the initial time and the boundary.
Ma and Wang \cite{MW} provided a  characterization of the boundary layer separation of 2-D incompressible
viscous fluids. They considered  a separation equation linking a separation
location and a time with the Reynolds number, the external forcing and the
initial velocity field.
In the Dirichlet boundary condition case, which corresponds to the
boundary layer separation, the above mentioned work of Ma and Wang provides detailed information on the  flow transition near
the critical time.
However, none of the above studies has shown the mechanism behind the reverse flow phenomena rigorously.
In this paper we observe a mathematical clue to analyze such reverse flow phenomena. More precisely, we consider the non-stationary Navier-Stokes equations with an initial datum having a parallel laminar  flow (we define rigorously later).
We show that the
direction of the material differentiation is opposite to the initial flow direction
and
 effect of the material differentiation (inducing the reverse flow) becomes bigger
when
the
curvature of the boundary becomes bigger.
Before showing such result,  we show the parallel laminar flow cannot be a stationary Navier-Stokes flow near a portion of the boundary with nonzero curvature.

Now we formulate our results in the following mathematical setting.
Let $\Omega$ be a either bounded or unbounded domain with smooth boundary in $\mathbb R^2$. The stationary Navier-Stokes
 equations are expressed as:
\begin{equation}\label{Euler and NS}
-\nu\Delta u+(u\cdot\nabla)u=-\nabla p\quad\text{and}\quad \div\ u=0\quad in \quad \Omega
\end{equation}
with $u|_{\partial\Omega}=0$.
We need to handle a shape of the boundary $\partial\Omega$ precisely, thus we set a parametrized smooth boundary $\varphi :(-\infty,\infty)\mapsto \mathbb R^2$ as
$
 |\partial_s\varphi(s)|=1,
$
$\theta(\partial_s\varphi(s))$ is a decreasing function,
$
\cup_{-\infty<s<\infty}\varphi(s)\subset\partial\Omega,
$
where $\theta(w)$ is defined by a vector $w=r(\cos\tilde\theta,\sin \tilde\theta)$,
$\theta(w):=\tilde \theta$.
For a technical sense, we need to assume that there are
$\bar s_0$ and $\bar s_2$ (we set $\bar s_1=0$, $\bar s_0<\bar s_1<\bar s_2$)
s.t.
$\varphi(\bar s_0)=(0,0)$, $\theta(\partial_s\varphi(s))|_{s=\bar s_0}=0$,
 $|\partial_s^2\varphi(s)|$ is monotone increasing for $s\in[\bar s_0, \bar s_1]$, $|\partial_s^2\varphi(s)|=1/\delta$ for $s\in[\bar s_1,\bar s_2]$, where $1/\delta$ is a constant curvature of a part of obstacle boundary $\cup_{s\in[\bar s_1,\bar s_2]}\varphi(s)$.
Note that there are $\tilde s$ and $\tilde x$, then $\varphi(s)$ ($s\in[\bar s_1,\bar s_2]$) can be expressed as
\begin{equation}\label{cone shape}
\varphi(s)=\delta\left(\cos\frac{(s+\tilde s)}{\delta},\sin\frac{(s+\tilde s)}{\delta}\right)+\tilde x.
\end{equation}
Here we mainly consider  laminar type flows (we define the ``laminar type flow"  later) in a localized region near the boundary, thus we need to   assume $u\not=0$ (no stagnation point) in
$\Omega$ near $\cup_{\bar s_0<s<\bar s_2}\varphi(s)$.
Boundary layers appear on the surface of bodies in viscous flow because the fluid seems to stick to the boundary $\partial\Omega$. Right at the surface the flow has zero relative speed and this fluid transfers momentum to adjacent layers through the action of viscosity. We need to express such phenomena in pure mathematics.
To do so,  we need to  give the following coordinates.
\begin{definition} (Normal coordinate.)
For $s\in [\bar s_0,\bar s_2]$, let
$$
\Phi(s,r)=\Phi_\varphi(s,r):=(\partial_s\varphi(s))^{\perp}r+\varphi(s).
$$
We define $\perp$ as the upward direction.
\end{definition}
\begin{remark}\label{Poiseuille}\ (An infinite channel with uniform width.)
We give a typical example of $\Omega$.
Assume that there is $\bar s_3(>\bar s_2)$ s.t. $|\partial_s^2\varphi(s)|$ is monotone decreasing for $s\in[\bar s_2,\bar s_3]$
and $|\partial_s^2\varphi(\bar s_0)|=|\partial_s^2\varphi(\bar s_3)|=0$.
Assume moreover that $\Omega$ is a straight pipe  except near the origin, namely, we assume that
$\varphi(s)=(s,0)$ for $s<0$,
$\theta(\partial_s\varphi(s))$ is a  constant for $s>\bar s_3$ and
$-\pi<\theta(\partial_s\varphi(\bar s_3))<0$.
Then we find an infinite channel with uniform width as
\begin{equation*}
\Omega=\{\Phi(s,r):s\in(-\infty,\infty),r\in(0,1)\}.
\end{equation*}
In this case, there are several results related to the Poiseuille flow (see \cite{K} for example).
\end{remark}

\begin{remark}\label{Phi}
By a scaling argument with a cone shape \eqref{cone shape}, we see that
for any  $s_1,s_2\in[\bar s_1,\bar s_2]$ with $s_1<s_2$,
\begin{equation*}
|\cup_{s'=s_1}^{s_2}\Phi(s',r_1)|=\left(\frac{r_{1}+\delta}{r_{2}+\delta}\right)|\cup_{s'=s_1}^{s_2}\Phi(s',r_2)|.
\end{equation*}
\end{remark}

In order to state main theorems, we need to give several definitions.

\begin{definition} (Normalized streamline and its direction.)
Let
$\gamma_a$ be in $\Omega$ near $\cup_{\bar s_0<s<\bar s_2}\varphi(s)$
which satisfies
$$
\partial_s\gamma_a(s)=\left(\frac{u}{|u|}\right)(\gamma_a(s)),\quad
\gamma_a(0)=a\in \Omega\quad\text{near}\quad \cup_{\bar s_0<s<\bar s_2}\varphi(s).
$$
Moreover, we assume that the flow is to the rightward direction (laminar flow direction), namely,
\begin{equation*}
0<\left\langle\partial_s\varphi(s),\frac{u}{|u|}(\Phi(s,r))\right\rangle<1
\end{equation*}
for $s\in [\bar s_0,\bar s_2]$ and sufficiently small $r$,
where $\langle\cdot,\cdot\rangle$ means inner product.
\end{definition}
We have already defined that there is no stagnation point near the boundary.
Thus the above definition is well defined.
\begin{definition} (Poincar\'e map.)\
For fixed $s$ and $s_1$ sufficiently close to each other,
let $s_{min}$ be the minimum of $s'>0$ for which there exists $\tau=\tau(s')$ such that
$\Phi(s_1,\tau(s'))=\gamma_{\Phi(s,r)}(s')$.
Let $L(r)=L_{s,s_1}(r)=\tau(s_{min})$.
\end{definition}
\begin{definition} (Classification of laminar type flows.)\
The velocity $u$ near the boundary $\cup_{s'\in[\bar s_0,\bar s_2]}\varphi(s')$ is called:

``Strong diverging laminar flow" iff $L(r)/r>C>1$ near the boundary,

``Weak diverging laminar flow" iff $L(r)/r\geq 1$ and $R(r)/r\to 1$ $(r\to 0)$,

``Parallel laminar flow" iff $L(r)/r=1$ near the  boundary.
\end{definition}
\begin{remark}\label{parallel case laminar flow}
In the parallel laminar flow case, we see that  $|u(\gamma_a(s))|=|u(a)|$ if $a$ and
$\gamma_a(s)$ are sufficiently close to $\cup_{s\in [\bar s_0,\bar s_2]}\varphi(s)$.
\end{remark}

In this paper we only consider the  parallel laminar flow case.
In the parallel laminar flow case, we see that
\begin{equation}\label{parallel direction}
\theta(\partial_s\gamma_{(0,x_2)}(s))|_{s=0}=0\quad\text{for sufficiently small}\quad  x_2>0.
\end{equation}
In this case, it is natural to  set
\begin{equation}\label{wake function}
u_1(0,x_2)=\alpha_1x_2-\frac{\alpha_2}{2}x_2^2(=:h(x_2))\quad\text{and}\quad u_2(0,x_2)=0\quad (\alpha_1,\alpha_2>0)
\end{equation}
in order to satisfy \eqref{Euler and NS} for sufficiently small  $x_2>0$. In this case, 
we can say that the boundary layer thickness is $\alpha_1/\alpha_2$.
This setting is based on  ``Poiseuille flow" (see Remark \ref{Poiseuille} and  \cite{Kob} for example).
Thus the parallel flow satisfying \eqref{wake function} is one of the candidate of the solution to the stationary Navier-Stokes equation.

\begin{remark}\label{uniform velocity}
We see from Remark \ref{parallel case laminar flow} that
\begin{equation}\label{equivalent to the initial position}
h(x_2)=|u(\Phi(s,x_2))|
\end{equation}
if $\Phi(s,x_2)$ is sufficiently close to $\cup_{s\in[\bar s_0, \bar s_2]}\varphi(s)$.
\end{remark}

The first main result is as follows:
\begin{theorem}\label{stationary}\ (Stationary case.)
Let us choose $\nu>0$ arbitrarily and fix it.
Then there is no   parallel laminar flow satisfying  the stationary Navier-Stokes equation \eqref{Euler and NS} and \eqref{wake function} near $\cup_{s'=\bar s_1}^{\bar s_2}\varphi(s)$.
\end{theorem}

Next we give a result in the non-stationary Navier-Stokes case.
The non-stationary Navier-Stokes equations are expressed as

\begin{equation}\label{NS}
\begin{cases}
\partial_t u-\nu\Delta u+(u\cdot\nabla)u=-\nabla p,\quad  \div\ u=0\quad in \quad \Omega\subset\mathbb{R}^2,\\
u|_{\partial \Omega}=0\quad\text{and}\quad u|_{t=0}=u_0.
\end{cases}
\end{equation}
Throughout the non-stationary case we assume that  $p(x,t)|_{t=0}$ is smooth in space valuable. If
$u_0$ is smooth and in $L^2(\Omega)$ (which means finite energy), then the solution $(u,p)$ is a global-in-time smooth solution (see \cite{L}).
By Theorem \ref{stationary}, we see that if we take an initial datum which has a parallel laminar flow near the boundary,
the solution $u$ must be a non-stationary flow in some small time interval near the initial time. In  this case,
we observe a mathematical clue to the separation phenomena. More precisely, we see that
the
direction of the material differentiation is opposite to the initial flow direction
and
the effect of the material differentiation
becomes bigger  
when 
the curvature of the boundary becomes bigger.
Let be more precise.
The initial datum having a parallel laminar flow satisfies
\begin{equation}\label{laminar}
\left\langle\frac{u_0}{|u_0|}(x),\partial_s\varphi(s)\right\rangle\to 1\quad\text{as}\quad x\to\varphi(s)\quad (s\in[\bar s_1,\bar s_2]).
\end{equation}
At some time $t>0$, if the solution $u(x,t)$ has the  reverse direction near the
boundary against the laminar flow direction,
we should have
\begin{equation}\label{reverse}
\left\langle\frac{u}{|u|}(x,t),\partial_s\varphi(s)\right\rangle\to -1\quad\text{as}\quad x\to\varphi(s)
\quad (s\in[\bar s_1,\bar s_2])\quad
\quad\text{for some}\quad t>0.
\end{equation}
To obtain such changing direction property,
we can expect that  the solution should satisfy
\begin{equation}\label{reverse differentiation}
\frac{\left\langle D_t u(x,t),\partial_s\varphi(s)\right\rangle}{\left\langle u_0(x),\partial_s\varphi(s)\right\rangle}\leq C<0\quad\text{as}\quad x\to\varphi(s)
\quad (s\in[\bar s_1,\bar s_2])
\end{equation}
for a long time interval up to the time $t$, where $D_tu:=\partial_t u+(u\cdot\nabla)u$ is the material differentiation.
We  show that the solution $u(x,t)$ to the non-stationary Navier-Stokes equation satisfies \eqref{reverse differentiation}
at the initial time, if the initial datum
has a parallel laminar flow near the boundary.
The following is the second main theorem.
\begin{theorem}
Assume that
$\nabla p(t,x)|_{t=0}$
is smooth.
If the initial datum satisfies \eqref{wake function} and has ``Parallel laminar flow", then we have
\begin{equation*}
\lim_{x\to \varphi(s)}\frac{\left\langle D_t u(x,t)|_{t=0},\partial_s\varphi(s)\right\rangle}{\langle u_0(x), \partial_s\varphi(s)\rangle}
= -\frac{\nu\alpha_2}{\delta\alpha_1}-\frac{\nu}{\delta^2}<0
\quad\text{for}\quad 
s\in [\bar s_1,\bar s_2].
\end{equation*}
\end{theorem}
Recall that the value $\alpha_1/\alpha_2$ is defined as the boundary layer thickness.
In general, the pressure is a non-local operator. However the values in the above theorem are depending only on the
behavior of the flow near the boundary  $\cup_{s\in[\bar s_1,\bar s_2]}\varphi(s)$.
\begin{remark}
There are direct and indirect evidences for the validity of the ``Kutta condition" in restricted regions
(see \cite{C}).
The method used in the above theorem may give another support for the validity of the Kutta condition  in pure mathematical sense.
\end{remark}

\begin{remark}
For the Poiseuille flow case,
the pressure is depending on $\nu$. More  precisely, if
\begin{equation*}
\Omega=\{x\in\mathbb{R}^2:0<x_2<1\},
\end{equation*}
then $u=(x_2-x_2^2,0)$, $p=(\nu/2)x_2$ is one of the stationary solution to \eqref{Euler and NS}.
To the contrary, for the whole space case (without any boundary effect), the pressure (in the sense of ``mild solution") is expressed as
\begin{equation*}
p=\sum_{i,j=1,2}R_iR_ju_iu_j,
\end{equation*}
where $R_j$ is the Riesz transform (see \cite{GIM,Ka, NY}).
In this observation, the dependence on the viscosity appears due to the Dirichlet boundary.
\end{remark}

\section{Proofs of the main theorems}
\noindent
{\it Proof of the first main theorem.}\ First we estimate $\Delta u$. For $\overline{s} \in [\overline{s}_{1} , \overline{s}_{2}]$ that is selected, we call the point $\varphi (\overline{s})$ to be $Q$. That is $Q = \varphi (\overline{s})$. Then, we  use the orthonormal frame $e_{1} = \partial_{s} \varphi |_{\overline{s}}$, $e_{2} = (\partial_{s} \varphi)^{\perp} |_{\overline{s}}$ at the point $Q$ to construct a cartesian coordinate system with the new $x_1$-axis to be the straight line which passes through $Q$ and parallel to the vector $e_{1}$, and the new $x_2$-axis to be the straight line which passes through $Q$ and is parallel to the vector $e_{2}$.
Then, for the given vector field $u$ in an open neighborhood near $\partial \Omega$, we define for each $y \in  \Omega$ near $\partial \Omega$, the two components of $u(y)$ with respect to the $e_{1}$ direction and the $e_{2}$ direction.
It is expressed as
\begin{equation*}
u(y)=v^1(y)e_1+v^2(y)e_2,
\end{equation*}
where $v^1(y):=\langle u(y),e_1\rangle$ and $v^2(y):=\langle u(y),e_2\rangle$.
Since $\varphi ([\overline{s}_{1}, \overline{s}_{2}])$ is known to be a circular arc with radius $\delta$, let $C$ to be the point at which the center of the circular arc $\varphi ([\overline{s}_{1}, \overline{s}_{2}])$ is located.
 Let $(s,r)$ to be the coordinate representation of the point $y$ in the coordinate system \emph{based at $Q$} which is specified by the orthonormal frame $\{e_{1}, e_{2}\}$.
That is, the point $y$ can be realized as $y = Q + s e_{1} + r e_{2}$. Since the vector field $u$ is a parallel laminar flow near $\partial \Omega$, the streamline of $u$ which passes through $y$ should be a circular arc with radius $|\overline{Cy}|$ and centered at $C$ also. Here, $|\overline{Cy}|$ is length of the line segment $\overline{Cy}$ That is, the distance between $C$ and $y$. Moreover, we  have
\begin{equation*}
|u(y)| = h (|\overline{Cy}| - \delta )
\end{equation*}
 with $\delta$ to be the radius of the circular arc $\varphi ([\overline{s}_{1}, \overline{s}_{2}])$, and $h(r) = \alpha_{1}r -\frac{\alpha_{2}}{2}r^{2}$. So, in order to compute $|u(y)|$, it is enough to compute $|\overline{Cy}|$ by means of geometry. Let $L_{x}$ to be the straight line which passes through $Q$ and is parallel to the $e_{1}$ direction. That is $L_{x}$ is the new $x-axis$ of the new coordinate system. First, let $A$ to be the unique point on $L_{x}$ such that the line segment $\overline{QA}$ is perpendicular to $\overline{Ay}$. Let $D$ to be the point of intersection between the line segment $\overline{Cy}$ and the line $L_{x}$. Observe that $|\overline{Ay}| = r$, $|\overline{QA}| = s$, and $|\overline{CQ}| = \delta$. We first compute $|\overline{QD}|$ and $|\overline{DA}|$ through the following observation. Since $\triangle DAy$ and $\triangle DQC$ are similar triangles, it follows that
\begin{equation*}
\frac{|\overline{DA}|}{r} = \frac{|\overline{QD}|}{\delta}.
\end{equation*}
Hence, by substituting $|\overline{DA}| = (r/\delta) |\overline{QD}|$ into the equation $s = |\overline{QD}| + |\overline{DA}|$, it follows that
\begin{equation*}
|\overline{DA}|  = \frac{rs}{\delta + r} \quad\text{and}\quad
  |\overline{QD}|  = \frac{\delta s}{\delta + r}.
\end{equation*}
This further gives
\begin{equation*}
|\overline{Dy}| = r \left\{ 1 + \left(\frac{s}{\delta + r}\right)^2 \right\}^{\frac{1}{2}} \quad\text{and}\quad
|\overline{DC}| = \delta \left\{ 1 + \left(\frac{s}{\delta + r}\right)^2 \right\}^{\frac{1}{2}}.
\end{equation*}
Hence, we have
\begin{equation*}
|\overline{Cy}| = |\overline{Dy}| + |\overline{DC}| = \{ (\delta + r)^{2} + s^{2}  \}^{\frac{1}{2}}.
\end{equation*}
As a result, we have $|u(p)| = h ( \{ (\delta + r)^{2} + s^{2}  \}^{\frac{1}{2}} - \delta  )$, and that
\begin{equation*}
\begin{split}
v^{1} &= h ( \{ (\delta + r)^{2} + s^{2}  \}^{\frac{1}{2}} - \delta  )\cdot \cos (\theta (y)) \\
v^{2} & = h ( \{ (\delta + r)^{2} + s^{2}  \}^{\frac{1}{2}} - \delta  )\cdot \sin (\theta (y)),
\end{split}
\end{equation*}
where, $\theta (y)$ is defined by the inscribed angle between $\overline{Dy}$ and $\overline{Ay}$. Since $\triangle DAy$ is a right angled triangle, it follows that
\begin{equation}\label{sincos}
\begin{split}
\sin (\theta (y)) &= - \sin (|\theta (y)|) =  -\frac{|\overline{DA}|}{|\overline{Dy}|} = -\frac{s}{\{ (\delta + r)^2 + s^2  \}^{\frac{1}{2}}} \\
\cos (\theta (y)) &= \cos (|\theta (y)|) = \frac{\delta + r }{\{ (\delta + r)^2 + s^2  \}^{\frac{1}{2}}} .
\end{split}
\end{equation}
Hence, it follows that
\begin{equation}\label{simple}
\begin{split}
v^{1}=v^1(Q+se_1+re_2) &= h ( \{ (\delta + r)^{2} + s^{2}  \}^{\frac{1}{2}} - \delta  )\cdot  \frac{\delta + r }{\{ (\delta + r)^2 + s^2  \}^{\frac{1}{2}}}\\
v^{2}=v^2(Q+se_1+re_2) & = - h ( \{ (\delta + r)^{2} + s^{2}  \}^{\frac{1}{2}} - \delta  )\cdot \frac{s}{\{ (\delta + r)^2 + s^2  \}^{\frac{1}{2}}}.
\end{split}
\end{equation}
Therefore we obtain, through direct computation that,
\begin{equation*}
\partial_r^2v^1=-\alpha_2,\ \partial_r^2v^2=0, \partial_{s}^{2} v^{1}|_{s= 0} = \frac{1}{r+\delta} \left(\alpha_{1} -\alpha_2 r- \frac{h(r)}{r+\delta}\right)
\quad\text{and}\quad
\partial_s^2v^{2}=0.
\end{equation*}
From \eqref{equivalent to the initial position}, we see $u(\gamma_{\Phi(0,r)}(s))=u(\Phi(s,r))=h(r)$. Thus the material differentiation can be calculated as
\begin{eqnarray*}
\left( (u\cdot \nabla)u \right)(\gamma_{\Phi(0,r)}(\bar s))&=&
\partial_s\left(u(\gamma_{\Phi(0,r)} (s))\right)|_{s=\bar s}\\
&=&
|u(\gamma_{\Phi(0,r)}(s))|\partial_s(\frac{u}{|u|}(\gamma_{\Phi(0,r)}(s)))|_{s=\bar s}\\
& &\ \ \ +
\partial_s|u(\gamma_{\Phi(0,r)}(s))|(\frac{u}{|u|}(\gamma_{\Phi(0,r)}(s)))|_{s=\bar s}\\
&=&
|u(\gamma_{\Phi(0,r)}(s))|\partial_s(\frac{u}{|u|}(\gamma_{\Phi(0,r)}(s)))|_{s=\bar s}\\
&=&
-\frac{1}{r+\delta} h(r)(\partial_s\varphi)^{\perp}(\bar s).
\end{eqnarray*}
Note that the above decomposition (into ``curvature part" and ``acceleration part") is already done in \cite{CCW}.
Therefore
\begin{equation}\label{wake side}
\nu\Delta u(\Phi(\bar s, r))-\left((u\cdot\nabla)u\right)(\Phi(\bar s, r))= P(r) \partial_{ s}\varphi(\bar s)
+P^\perp(r)(\partial_{ s}\varphi)^{\perp}(\bar s)
\end{equation}
for $\bar s\in[\bar s_1,\bar s_2]$, where
\begin{equation}\label{P}
P(r)=
  \nu\left (\frac{\alpha_{1}}{r+\delta} - \frac{\alpha_2 r}{r+\delta}-\frac{h(r)}{(r+\delta)^{2}} - \alpha_{2}\right),\quad
P^\perp(r):=\frac{h(r)}{r+\delta}.
\end{equation}
In order to see the pressure precisely,  we need the following ``normalized pressure-line" and ``boundary of the  level set of the pressure".
Let
\begin{equation*}
\partial_sq_a(s)=\frac{\nabla p}{|\nabla p|}(q_a(s)),\quad \partial_r q^{\perp}_a(r)=\left(\frac{\nabla p}{|\nabla p|}\right)^{\perp}(q^{\perp}_a(r))\quad\text{for}\quad|\nabla p|\not=0
\end{equation*}
with $q_a(0)=a$ and $q^{\perp}_a(0)=a$.
Moreover let us define a re-parametrized $q$ as follows:
\begin{definition}
(Poincar\'e map on the pressure lines.)\
For sufficiently small $\epsilon>0$, let $s_{min}$ be the minimum of $s'>0$ for which there exists $\tau=\tau(s')$ such that
\begin{equation*}
q_{\Phi(s,r)}(s')=q^{\perp}_{\Phi(s+\epsilon,r)}(\tau(s')).
\end{equation*}
Let $\eta(\epsilon)=\eta(\epsilon,s,r)=q_{\Phi(s,r)}(s_{min})$.
We denote $\cup^\epsilon\eta(\epsilon')$ as $\cup_{\epsilon'=0}^\epsilon\eta(\epsilon',s,r)$.
\end{definition}

Now we consider  the level set of $p$ precisely.
We choose the following finite points:
 $\{ \eta(\epsilon'(0)),\eta(\epsilon'(1)), \eta(\epsilon'(2)),\cdots,\eta(\epsilon'(N-1)), \eta(\epsilon'(N)) \}$ in order to satisfy that
$\epsilon'(0)=0$, $\epsilon'(N)=\epsilon$ and
 $|\eta(\epsilon'(\ell))-\eta(\epsilon'(\ell-1))|$ $(\ell=1,\cdots N)$ are the same value.
Then we set $\nabla p$ (in the discrete setting) from the level set of $p$
\begin{equation*}
\nabla p(\eta( \epsilon'(\ell)))\approx
 \frac{ p(\eta(\epsilon'(\ell))) - p(\eta(\epsilon'(\ell+1)))}
{|\eta(\epsilon'(\ell)) -\eta(\epsilon'(\ell+1))|}
 \frac{\eta(\epsilon'(\ell)) -\eta(\epsilon'(\ell+1))}
{|\eta(\epsilon'(\ell)) -\eta(\epsilon'(\ell+1))|}
\end{equation*}
and let $\vec t$ be a tangent vector defined as
\begin{equation*}
\vec t(\epsilon'(\ell))=
 \frac{\eta(\epsilon'(\ell)) -\eta(\epsilon'(\ell+1))}
{|\eta(\epsilon'(\ell)) -\eta(\epsilon'(\ell+1))|}.
\end{equation*}
Then we have in the discrete setting (piecewise linear approximation):
\begin{eqnarray*}
\int_{\cup^\epsilon\eta(\epsilon')}\nabla p\cdot \vec  t&= &\lim_{N\to \infty}\sum_{\ell=1}^N
\nabla p(\epsilon'(\ell))\cdot \vec  t(\epsilon'(\ell))|\eta(\epsilon'(\ell)) -\eta(\epsilon'(\ell+1))|\\
&=&
p(\eta(\epsilon,s, r))-p(\Phi( s, r)).
\end{eqnarray*}
Therefore
\begin{eqnarray}\label{widetilde p}
|\nabla p(\Phi( s, r))|&=&\lim_{\epsilon\to 0}\frac{1}{|\cup^\epsilon\eta(\epsilon')|}
\left|\int_{\cup^\epsilon\eta(\epsilon')}\nabla p\cdot \vec t\right|\\
\nonumber
&=&
\lim_{\epsilon\to 0}\frac{1}{|\cup^\epsilon\eta(\epsilon')|}\left|p(\eta(\epsilon, s, r))-p(\Phi( s, r))\right|.\\
\nonumber
\end{eqnarray}
Since  the direction of $\nabla p$, namely $\frac{\nabla p}{|\nabla p|}$ is determined by
$\frac{\Delta u}{|\Delta u|}$ and $\frac{\partial_s u(\gamma(s))}{|\partial_s u(\gamma(s))|}$,
and $u$ has a parallel laminar flow,  we see that
\begin{equation*}
\left\langle\frac{\nabla p}{|\nabla p|}\left(\Phi(s,r)\right),\frac{\partial_s\Phi(s,r)}{|\partial_s\Phi(s,r)|}\right\rangle
\end{equation*}
is independent of $s$. Thus there are  positive  $\hat s(r)$  and $\hat r(r)$ (independent of $s$) such that
\begin{equation*}
p(\Phi(s,r))=p(q^{\perp}_{\Phi(s,r)}(0))=p(q^\perp_{\Phi(s,r)}(\hat r(r)))=p(\Phi(s-\hat s(r),0))
\end{equation*}
and
\begin{equation*}
p(\eta(\epsilon,s,r))=p(q^{\perp}_{\Phi(s+\epsilon,r)}(0))=p(q^\perp_{\Phi(s+\epsilon,r)}(\hat r(r)))=p(\Phi(s+\epsilon-\hat s(r),0)).
\end{equation*}
Since
\begin{equation*}
\nabla p (\Phi(s,0)) = \nu \triangle u (\Phi (s, 0)) = \nu (\frac{\alpha_{1}}{\delta} - \alpha_{2}) \partial_{s} \varphi
\end{equation*}
derived from the boundary condition, we see
\begin{eqnarray*}
|p(\Phi(s+\epsilon-\hat s(r),0))-p(\Phi(s-\hat s(r),0))|&=&\left|\int_0^\epsilon\nabla p(\Phi(s',0))\cdot \varphi(s')\ ds'\right|\\
&=&
 \epsilon \nu |\frac{\alpha_{1}}{\delta} - \alpha_{2}|.
\end{eqnarray*}
Thus we have the following pressure estimate:
\begin{equation*}
|\nabla p(\Phi(s,r))|=
\lim_{\epsilon\to 0}\frac{\epsilon\nu|\alpha_1/\delta-\alpha_2|}{|\cup^\epsilon\eta(\epsilon')|}.
\nonumber
\end{equation*}
On the other hand,
since  $\theta\left((\partial_s\Phi)( s+\epsilon)\right)>\theta\left(\partial_s q_{\Phi(s,r)}(s+\epsilon)\right)$ for any $\epsilon$,
there is a negative number
 $\tau=\tau(\epsilon,s,r)<0$ such that
$q^{\perp}_{\Phi(s+\epsilon, r)}(\tau)=\eta(\epsilon, s, r)$.
Thus $D_\epsilon:=\cup_{\epsilon'=0}^\epsilon\cup_{r'=\tau(\epsilon',s,r)}^0
q^\perp_{\Phi( s+\epsilon', r)}(r')$ is well defined and
$\lim_{\epsilon\to 0}D_\epsilon/|D_\epsilon|$ is a right triangle
since
$\langle(\partial_s\Phi)( s, r),((\partial_\epsilon\eta)(0,s, r))\rangle$ does not change for any $\epsilon>0$
and\\
$\langle\partial_{\tau}q^\perp_{\Phi( s+\epsilon,r)}(\tau),(\partial_\epsilon\eta)(\epsilon, s, r)\rangle=0$ for any $\epsilon>0$.
It means that angles of the three corners  do not change for any $\epsilon$
and one of them is $\pi/2$.
More precisely,  among the three sides of the triangular region $D_{\epsilon}$, the side
$\cup_{0 \leq \epsilon' \leq \epsilon} \eta (\epsilon' , s, r )$, which is a part of the pressure line passing through $\Phi (s,r)$, is \emph{perpendicular to} the level set
of the function $p$, which passes through $\Phi (s + \epsilon , r)$, \emph{at the point of intersection $\eta (\epsilon , s , r )$}. In other words, we have the following observation.

\begin{itemize}
\item the right angle $\frac{\pi}{2}$ of the triangle $D_{\epsilon}$ is located at  the vertex $\eta (\epsilon , s ,r)$, which is the point of intersection between the side $\cup_{0 \leq \epsilon' \leq \epsilon} \eta (\epsilon' ,s,r )$ and the level set of $p$ passing through $\Phi(s+ \epsilon ,r)$.
\end{itemize}

 So, the arc segment
$|\cup_{s \leq s' \leq s + \epsilon} \Phi (s' , r)|$ should be the \emph{longest side} of the triangle $D_{\epsilon}$ whose \emph{opposite angle} in $D_{\epsilon}$ is the \emph{right angle} $\frac{\pi}{2}$ located at the vertex $\eta (\epsilon , s, r)$ of $D_{\epsilon}$. According to the above reasoning, it is not hard to see that
\begin{equation}\label{correct}
\frac{|P(r)|}{[(P(r))^{2} + (P^{\perp}(r))^{2}]^{\frac{1}{2}}} = \lim_{\epsilon \rightarrow 0^{+}} \frac{|\cup_{0 \leq \epsilon' \leq \epsilon} \eta (\epsilon' ,s,r ) |}{|\cup_{s \leq s' \leq s + \epsilon} \Phi (s' , r) |} .
\end{equation}
In accordance with the relation in \eqref{correct}, it follows that
\begin{equation*}
\begin{split}
\lim_{\epsilon \rightarrow 0^{+}} \frac{\epsilon}{ |\cup_{0 \leq \epsilon' \leq \epsilon} \eta (\epsilon' ,s,r ) |  } & =
\lim_{\epsilon \rightarrow 0^{+}} \frac{\epsilon}{ |\cup_{s \leq s' \leq s + \epsilon} \Phi (s' , r)|   } \cdot \frac{|\cup_{s \leq s' \leq s + \epsilon} \Phi (s' , r)| }{|\cup_{0 \leq \epsilon' \leq \epsilon} \eta (\epsilon' ,s,r) |} \\
& = \frac{\epsilon}{(\frac{\delta + r}{\delta } |\cup_{s \leq s' \leq s + \epsilon} \Phi (s' , 0)|   )} \cdot \frac{[(P(r))^{2} + (P^{\perp}(r))^{2}]^{\frac{1}{2}}}{|P(r)|} \\
& = \left(\frac{\epsilon}{(\frac{\delta + r}{\delta}) \epsilon }\right) \frac{[(P(r))^{2} + (P^{\perp}(r))^{2}]^{\frac{1}{2}}}{|P(r)|} \\
& = \left(\frac{\delta}{\delta + r}\right) \frac{[(P(r))^{2} + (P^{\perp}(r))^{2}]^{\frac{1}{2}}}{|P(r)|}.
\end{split}
\end{equation*}
Hence, it follows from relation \eqref{widetilde p} and the above that
\begin{equation*}
|\nabla p (\Phi (s,r) )| = \nu |\frac{\alpha_{1}}{\delta} - \alpha_{2}| \left(\frac{\delta}{\delta + r}\right) \frac{[(P(r))^{2} + (P^{\perp}(r))^{2}]^{\frac{1}{2}}}{|P(r)|}.
\end{equation*}
However, since we naturally have $|\nabla p (\Phi (s,r) )| =  [(P(r))^{2} + (P^{\perp}(r))^{2}]^{\frac{1}{2}}$, it follows from the above relation that
\begin{equation*}
\nu |\frac{\alpha_{1}}{\delta} - \alpha_{2}|\left (\frac{\delta}{\delta + r}\right) \frac{[(P(r))^{2} + (P^{\perp}(r))^{2}]^{\frac{1}{2}}}{|P(r)|} = [(P(r))^{2} + (P^{\perp}(r))^{2}]^{\frac{1}{2}}.
\end{equation*}
The above expression, together with the expression $P(r) = \nu (\frac{\alpha_{1}}{r+\delta}- \frac{\alpha_2r}{(r+\delta)} - \frac{h(r)}{(r+\delta)^{2}} - \alpha_{2} )$ directly implies that
\begin{equation}\label{absurd}
 |\frac{\alpha_{1}}{\delta} - \alpha_{2}| \left(\frac{\delta}{\delta + r}\right) = \left|\frac{\alpha_{1}}{r+\delta}- \frac{\alpha_2r}{(r+\delta)} - \frac{h(r)}{(r+\delta)^{2}} - \alpha_{2} \right|.
\end{equation}
Our task here is to derive a contradiction from the relation in \eqref{absurd}. First, let us deal with the case in which $|\frac{\alpha_{1}}{\delta} - \alpha_{2}| > 0 $. In this case, since $(\frac{\alpha_{1}}{r+\delta} - \frac{\alpha_2r}{(r+\delta)} -\frac{h(r)}{(r+\delta)^{2}} - \alpha_{2} ) \rightarrow \frac{\alpha_{1}}{\delta} - \alpha_{2}$, as $r \rightarrow 0^{+}$, it follows that $(\frac{\alpha_{1}}{r+\delta} - \frac{\alpha_2r}{(r+\delta)}- \frac{h(r)}{(r+\delta)^{2}} - \alpha_{2} )$ will have the same sign as that of $\frac{\alpha_{1}}{\delta} - \alpha_{2}$, for all sufficiently small $r > 0$. That is, for all $r > 0$ to be sufficiently small, the following assertion holds.\\

 $(\frac{\alpha_{1}}{r+\delta}- \frac{\alpha_2r}{(r+\delta)} - \frac{h(r)}{(r+\delta)^{2}} - \alpha_{2} )$ is positive if and only if $\frac{\alpha_{1}}{\delta} - \alpha_{2}$ is positive.\\

\noindent
As a result, the relation in \eqref{absurd} can be rephrased as
\begin{equation}\label{absurd2}
\left (\frac{\alpha_{1}}{\delta} - \alpha_{2}\right) \left(\frac{\delta}{\delta + r}\right) = \frac{\alpha_{1}}{r+\delta}- \frac{\alpha_2r}{(r+\delta)} - \frac{h(r)}{(r+\delta)^{2}} - \alpha_{2} .
\end{equation}
However, by some algebraic computation, relation as given in \eqref{absurd2} is equivalent to saying that the following relation holds for all $r > 0$ sufficiently close to $0$,
\begin{equation*}
\alpha_{2} \left(\frac{\delta}{\delta + r} - 1\right) = \frac{h(r)}{(\delta + r)^2}+ \frac{\alpha_2r}{r+\delta}.
\end{equation*}
However, the above relation says that,  for $r > 0$ to be sufficiently small, the left hand side $\alpha_{2} [\frac{\delta}{\delta + r} - 1]$,which is strictly negative, \emph{equals} the right hand side $\frac{h(r)}{\delta + r} + \frac{\alpha_2r}{(r+\delta)}$, which is strictly positive. This is \emph{absurd}, and hence we have arrive at a contradiction in this case.\\
The case of $\frac{\alpha_{1}}{\delta} - \alpha_{2} = 0$ is even easier to handle, since in this case, relation \eqref{absurd} immediately gives
\begin{equation*}
0 = \frac{\alpha_{1}}{\delta + r}- \frac{\alpha_2r}{r+\delta} - \frac{h(r)}{(\delta + r)^2} - \alpha_{2} =   \frac{\delta \alpha_{2}}{\delta + r} - \frac{h(r)}{(\delta + r)^2} - \alpha_{2} 
\end{equation*}
which eventually leads to the relation $ \alpha_{2} [\frac{\delta}{\delta + r} - 1] = \frac{h(r)}{(\delta + r)^2}+ \frac{\alpha_2r}{(r+\delta)} $, which is again absurd. So, we also have a contradiction.

\vspace{0.5cm}

\noindent
{\it Proof of the second main theorem (Non-stationary case.)}\
Let $p(x):=p(x,t)|_{t=0}$ and $u(x):=u(x,t)|_{t=0}$.
From \eqref{P}, we see that
\begin{equation*}
\Delta u(\Phi(s,r))=\nu\left(\frac{\alpha_1}{r+\delta}- \frac{\alpha_2r}{r+\delta}-\frac{h(r)}{(r+\delta)^2}-\alpha_2\right)\partial_s\varphi(s).
\end{equation*}
Thus we have
\begin{equation}\label{pressure boundary}
\nabla p(\Phi(s,r))|_{r=0}= \nu\Delta u(\Phi(s,r))|_{r=0}=\nu\left(\frac{\alpha_1}{\delta}-\alpha_2\right)\partial_s \varphi(s)
\end{equation}
derived from the Dirichlet boundary condition.
In order to estimate $\nabla p$ only using smoothness and \eqref{pressure boundary},
we use ``normalized pressure-line" and ``boundary of the level set of the pressure" defined in the stationary case.
The key is to estimate $|\cup_{\epsilon'=0}^\epsilon\zeta(\epsilon')|$ (we define it  later) which is similar to $|\cup_{\epsilon'=0}^\epsilon\eta(\epsilon')|$ defined in the stationary case.
\begin{definition}\  Here we give three definitions which are needed in a geometric observation.

\begin{itemize}

\item

Let $\hat s=\hat s(s,r)$ and $\hat r=\hat r(s,r)$ be such that
\begin{equation*}
q^\perp_{\varphi(s)}(\hat r)=\Phi(\hat s,r) \quad\text{for}\quad s\in[\bar s_1,\bar s_2].
\end{equation*}

\item

For sufficiently small $\epsilon>0$ 
,
 let $s_{min}$ be the minimum of $s'>0$ for which there exists
$\tau=\tau(s')$ such that $q_{F}(s')=q^{\perp}_{\varphi(s+\epsilon)}(\tau(s'))$, where
$F=q^{\perp}_{\varphi(s)}(\hat r(s,r))$.
Let $\zeta(\epsilon)=\zeta(\epsilon, s,r)=q_{F}(s_{min})$.
We denote $\cup^\epsilon\zeta(\epsilon')$ as $\cup_{\epsilon'=0}^\epsilon\zeta(\epsilon',s,r)$.

\item

Let $\hat{\hat s}=\hat{\hat s}(\epsilon,s,r)$ and $\hat{\hat{ r}}=\hat {\hat r}(\epsilon,s,r)$ be such that $\zeta(\epsilon, s,r)=\Phi(\hat{\hat s}, \hat{\hat r})$. 
\end{itemize}
\end{definition}

\vspace{0.5cm}

Due to the smoothness of $p$ and \eqref{pressure boundary}, we can estimate $\hat s$ and $\hat{\hat s}$ as
\begin{equation}\label{error estimate}
s-cr^2\leq \hat s\leq s+c r^2
\quad\text{and}\quad
s+\epsilon-c_1\hat{\hat r}^2\leq \hat {\hat s}\leq s+\epsilon+c_2\hat{\hat r},
\end{equation}
where $c$ is a positive constant independent of $\epsilon$, $s$, $r$ and $\hat {\hat r}$.
Let us set $\hat D_{\epsilon,r}$ (which is similar to $D_\epsilon$ in the stationary case but different) as follows:
\begin{equation*}
\hat D_{\epsilon,r}:=\cup_{\epsilon'=0}^\epsilon\cup_{r'=0}^r\zeta(\epsilon',s,r').
\end{equation*}
Since all of the four angles are $\pi/2$,
\begin{equation*}
\lim_{\stackrel{\epsilon,r\to 0}{r/\epsilon=const.}}\frac{\hat D_{\epsilon,r}}{|\hat D_{\epsilon,r}|}
\end{equation*}
is a rectangular. Thus we can see that for sufficiently small $\hat \epsilon>0$, there are $R$, $\hat R$ and $\hat{\hat R}$ such that
\begin{equation}\label{rectangular}
(1-\hat\epsilon)r\leq \hat{\hat r}\leq (1+\hat \epsilon)r
\end{equation}
for $r<R$ with $r/\epsilon=const.$, $\hat r<\hat R$ and $\hat{\hat r}<\hat {\hat R}$.
From \eqref{error estimate}, we have 
\begin{equation*}
s+\epsilon-c(1+\hat \epsilon)^2r^2\leq \hat{\hat s}\leq s+\epsilon+c(1+\hat \epsilon)^2 r^2
\end{equation*}
and then
\begin{equation*}
\epsilon-cr^2-c(1+\hat \epsilon)r^2\leq \hat {\hat s}-\hat s\leq \epsilon+cr^2+c(1+\hat \epsilon)r^2.
\end{equation*}
Due to the smoothness of $p$ and \eqref{pressure boundary},
we also see that
\begin{equation*}
 1-cr^2\leq \left\langle\left(\frac{\nabla p}{|\nabla p|}\right)(\Phi(s,r))  ,
\partial_s   \varphi(s)\right\rangle
\leq 1.
\end{equation*}
By using a piecewise linear approximation,  we  can estimate   $|\cup^\epsilon\zeta(\epsilon')|$ as
\begin{eqnarray*}
|\cup^\epsilon\zeta(\epsilon')|
&=&
\lim_{N\to \infty}\sum_{n=0}^N\frac{\sum^n_{k=0} \Delta r_k+r}{r \cos\theta(\Phi(\hat s+n\Delta s,r+\sum_{k=0}^n\Delta r_k))}\Delta s,
\end{eqnarray*}
where $\Delta r_0=0$, $\Delta r_k=\tan \theta(\Phi(\hat s+(k-1)\Delta s,r+\sum^{k-1}_{k'}\Delta r_{k'}))$ for $k\geq 1$,
$\theta (\Phi(s,r))$ is defined as the angle between $(\nabla p/|\nabla p|)(\Phi(s,r))$ and $\partial_s\varphi(s)$,
and
we choose $\Delta s$  in order to satisfy
\begin{equation*}
N\Delta s=|\cup_{s'=\hat s}^{\hat{\hat s}}\Phi(s',r)|=\frac{\delta+r}{\delta}(\hat{\hat s}-\hat s).
\end{equation*}
Note that  $\Delta s$ tends to $0$ as $N\to \infty$. It means that $\Delta r_k$ tends to $0$ as $N\to\infty$ for fixed $k$.
Then we see that 
\begin{equation*}
\lim_{n\to\infty,  N\to \infty}\sum_{k=0}^n\Delta r_k+r=\hat {\hat r}.
\end{equation*} 
By the above estimates, we have the following lower and upper bound of $|\cup^\epsilon \zeta(\epsilon')|$: 
\begin{eqnarray}\label{lower}
|\cup^\epsilon\zeta(\epsilon')|
&\geq&
(1-\hat \epsilon)\frac{r+\delta}{\delta}\left(\epsilon-cr^2-c(1+\hat \epsilon)^2r^2\right),\\
|\cup^\epsilon\zeta(\epsilon')|
&\leq&
\frac{1+\hat \epsilon}{1-cr^2}\cdot\frac{r+\delta}{\delta}\left(\epsilon+cr^2+c(1+\hat\epsilon)^2r^2\right)
\nonumber
\end{eqnarray}
for sufficiently small $\epsilon$ and $r$ with $\epsilon/r=const$. Note that we can take $\hat \epsilon>0$ arbitrarily small. 
We have from the explicit representation of the gradient of the pressure (here we consider the case  $\alpha_1/\delta-\alpha_2>0$):
\begin{eqnarray*}
\langle \nabla p(\varphi(s)), \partial_s\varphi(s)\rangle
&=& |\nabla p(\varphi(s))|\left\langle\frac{\nabla p}{|\nabla p|}(\varphi(s)),\varphi(s)\right\rangle\\
&= &
|\nabla p(\Phi(\hat s,r))|\\
&=&\lim_{\epsilon\to 0}\frac{1}{|\cup^\epsilon \zeta(\epsilon')|}\left|\int_{\cup^\epsilon\zeta(\epsilon')}\nabla p(x)\cdot \vec t\ dx\right|\\
&=&
\lim_{\epsilon\to 0}\frac{1}{|\cup^\epsilon \zeta(\epsilon')|}\left|p(\varphi(s+\epsilon))-p(\varphi(s))\right|\\
&=&
\lim_{\epsilon\to 0}\frac{\epsilon\nu \left(\alpha_1/\delta-\alpha_2\right)}{|\cup^\epsilon \zeta(\epsilon')|}.
\end{eqnarray*}
Even in the case  $\alpha_1/\delta-\alpha_2\leq 0$, we have the same estimate as above.
Recall that
\begin{equation*}
\langle \nu\Delta u, \partial\varphi(s)\rangle=\nu\left(\frac{\alpha_1}{r+\delta} - \frac{\alpha_2r}{(r+\delta)}-\frac{h(r)}{(r+\delta)^2}-\alpha_2\right)
\end{equation*}
and
\begin{equation*}
\langle u_0,\partial_s\varphi(s)\rangle=\alpha_1r-\frac{\alpha_2}{2}r^2.
\end{equation*}
We already have the lower and upper bounds of $|\cup^\epsilon\zeta(\epsilon')|$  (see \eqref{lower}), we have that 
\begin{eqnarray*}
\frac{\langle D_tu,\partial_s\varphi(s)\rangle}{\langle  u_0,\partial_s\varphi(s)\rangle}&=&
\frac{\langle(\nu\Delta u-\nabla p),\partial_s\varphi(s)\rangle}{\langle u_0,\partial_s\varphi(s)\rangle}\\
&=&
\lim_{\stackrel{\epsilon,r\to 0}{\epsilon/r=const.}}\frac{\langle\nu\Delta u(\Phi(s,r)),\partial_s\varphi(s)\rangle-\left(\epsilon\nu(\alpha_1/\delta-\alpha_2)/|\cup^\epsilon\zeta(\epsilon')|\right)}{\langle u_0(\Phi(s,r)),\partial_s\varphi(s)\rangle }\\
&= &-\frac{\nu\alpha_2}{\delta\alpha_1}-\frac{\nu}{\delta^2}
\end{eqnarray*}
for  $s\in[\bar s_1,\bar s_2]$.
Therefore we have the   the desired estimate.

\vspace{0.5cm}

{\bf Acknowledgments.}\
The author would like to thank Professor Chi hin Chan for  giving me many advices.
The author would also like to thank Professor Takashi Sakajo for letting me know
the  article \cite{C}. The author is partially supported by JST CREST.

\bibliographystyle{amsart}

\end{document}